
\input amstex
\magnification=1200
\documentstyle{amsppt}
\leftheadtext{}
\rightheadtext{}
\NoBlackBoxes
\define\End{\operatorname{End}}
\define\Lie{\operatorname{Lie}}
\define\Hom{\operatorname{Hom}}
\define\card{\operatorname{card}}
\define\ad{\operatorname{ad}}
\redefine\P{\hat P_t}
\define\Q{\hat Q_t}
\define\R{\hat R_t}
\define\A{\hat A_t}
\redefine\B{\hat B_t}
\define\C{\hat C_t}
\define\spn{\operatorname{span}}
\define\vz{\text{\bf !`}}
\topmatter
\title Classical and quantum dynamics of noncanonically coupled oscillators
and Lie superalgebras
\endtitle
\author Denis V. Juriev
\endauthor
\affil\eightpoint Mathematical Division, Research Institute for System
Studies,\linebreak
Russian Academy of Sciences, Moscow, Russia\footnote{ Current address:
Laboratoire de Physique Th\'eorique de l'\'Ecole Normale Sup\'erieure, 24 rue
Lhomond, 75231 Paris Cedex 05, France (e-mail:
juriev\@physique.ens.fr).\newline}
\endaffil
\date funct-an/9409003 \enddate
\abstract The classical and quantum dynamics of noncanonically coupled
oscillators is investigated in its relation to Lie superalgebras. It is shown
that the quantum dynamics admits a hidden (super)hamiltonian formulation and,
hence, preserves the initial operator relations.
\endabstract
\endtopmatter
\document

\head I. INTRODUCTION
\endhead

It is well--known that classical and quantum dynamics of hamiltonian systems
is often described by remarkable algebraic structures such as Lie algebras,
their nonlinear generalizations and (quantum) deformations [1,2]. It seems that
not less important objects govern a behaviour of the interacting hamiltonian
systems and that they maybe unravelled in a certain way. There exist several
forms of an interaction of hamiltonian systems: often it has a potential
character, sometimes it is ruled by a deformation of the Poisson brackets;
however, one of the most intriguing and mathematically less explored forms is
a nonhamiltonian (noncanonical) interaction, which can not be described by
deformations of the standard hamiltonian data (Poisson brackets and
hamiltonians). Sometimes, such curious interaction is realized by the
dependence (which is linear in the simplest cases [3] being nonlinear in
general [4]) of the Poisson brackets of one hamiltonian system on the state of
another [3]. The pair of noncanonically coupled oscillators is one of the
simplest and the most crucial examples of the nonhamiltonian interaction [3].
The purpose of this paper is to describe the classical and quantum dynamics of
noncanonically coupled oscillators in a general setting. For that we use the
following algebraic structures (and their representations): (a) isotopic pairs
(a particular linear case of general I--pairs of the note [4]), (b)
anti--Jordan pairs, (c) anti--Lie triple systems, (d) Lie superalgebras. A
brief description of relations between these algebraic structures is presented
at the end of par.II.

\head II. GENERAL DEFINITIONS
\endhead

\definition{Definition 1 {\rm [3]}} The pair $(V_1,V_2)$ of linear spaces is
called {\it an isotopic pair\/} iff there are defined two mappings
$m_1:V_2\otimes\bigwedge^2V_1\mapsto V_1$ and
$m_2:V_1\otimes\bigwedge^2V_2\mapsto V_2$ such that the mappings $(X,Y)\mapsto
[X,Y]_A=m_1(A,X,Y)$ ($X,Y\in V_1$, $A\in V_2$) and $(A,B)\mapsto
[A,B]_X=m_2(X,A,B)$ ($A,B\in V_2$, $X\in V_1$) obey the Jacobi identity for
all values of a subscript parameter (such operations will be called {\it
isocommutators\/} and the subscript parameters will be called {\it isotopic
elements\/} or shortly {\it isotopies\/}) and are compatible to each other,
i.e. the identities
$$\align
[X,Y]_{[A,B]_Z}=&\tfrac12([[X,Z]_A,Y]_B+[[X,Y]_A,Z]_B+[[Z,Y]_A,X]_B-\\
&[[X,Z]_B,Y]_A-[[X,Y]_B,Z]_A-[[Z,Y]_B,X]_A)\endalign$$
and
$$\align
[A,B]_{[X,Y]_C}=&\tfrac12([[A,C]_X,B]_Y+[[A,B]_X,C]_Y+[[C,B]_X,A]_Y-\\
&[[A,C]_Y,B]_X-[[A,B]_Y,C]_X-[[C,B]_Y,A]_X)\endalign$$
($X,Y,Z\in V_1$,
$A,B,C\in V_2$) hold.
\enddefinition

Let's discuss this defintion.

First, it may be considered as a result of an axiomatization of the following
trivial construction: let $\Cal A$ be an associative algebra (f.e. any matrix
one) and $V_1$, $V_2$ be two linear subspaces in it such that $V_1$ is closed
under the isocommutators $(X,Y)\mapsto [X,Y]_A=XAY-YAX$ with isotopic elements
$A$ from $V_2$, whereas $V_2$ is closed under the isocommutators $(A,B)\mapsto
[A,B]_X=AXB-BXA$ with isotopic elements $X$ from $V_1$.

\remark{Remark 1} Let $H$ be a (finite dimensional) linear space. If $\Cal A$
is
a subspace of $\End(H)$ let's put $\Cal A^{{\vz}}=\{X\in\End(H),\forall
A\in\Cal
A, \forall B\in\Cal A, AXB-BXA\in\Cal A\}$. Then $\Cal A\subseteq\Cal
A^{{\vz}{\vz}}$ and $(\Cal A,\Cal A^{{\vz}})$ is an isotopic pair.
\endremark

\define\abZ{\alpha\underset Z\to\lozenge\beta}
It is rather interesting to unravel the most general setting for such
construction. Namely, let $[\cdot,\cdot]_\alpha$ ($\alpha\in\frak A$) be a
linear family of compatible Lie brackets on $V$. When the bracket
$[\cdot,\cdot]_{\abZ}$ ($Z\in V$) defined on $V$ as
$$[X,Y]_{\abZ}=\tfrac12([[X,Z]_A,Y]_B+[[X,Y]_A,Z]_B+[[Z,Y]_A,X]_B-
(\alpha\longleftrightarrow\beta))$$
is a Lie bracket compatible with brackets $[\cdot,\cdot]_\alpha$ and
$[\cdot,\cdot]_\beta$ ($\alpha, \beta\in\frak A$)?

Certainly, the $\frak g$--equivariant case is the most important one but,
unfortunately, I do not know any answer on this simple question. The existence
of $\lozenge$--operation on the space of compatible Lie brackets depends on
the correctness of the following {\it conjecture\/}: let
$[\cdot,\cdot]_\alpha$ ($\alpha\in\frak A$) be a linear family of compatible
Lie brackets on the linear space $V$, then there exist a linear space $H$ and
two linear mappings $T\in\Hom(V;\End(H))$ and $Q\in\Hom(\frak A;\End(H))$ such
that $T([X,Y]_\alpha)=T(X)Q(\alpha)T(Y)- T(Y)Q(\alpha)T(X)$.

\remark{Remark 2} Let $\Cal A\subseteq\End(H)$, $\dim\Cal A=n$, $\Lie_n$ be the
space of all Lie algebras of dimension $n$. Then there exists a natural
mapping $\Cal L:\Cal A^{{\vz}}\mapsto\Lie_n$. It should be mentioned that
$\card\Cal L(\Cal A^{{\vz}})$ may be not equal to $1$ so $\Cal L(X)$ and $\Cal
L(Y)$ are not the same in general for different $X$ and $Y$. It means that the
isocommutators $[\cdot,\cdot]_X$ and $[\cdot,\cdot]_Y$ determine structures of
nonisomorphic Lie algebras on the space $\Cal A$ in general (though they may
be isomorphic in particular).
\endremark

Second, one may compare def.1 with the definition of "anti--Jordan pairs" [5].
Namely,

\definition{Definition 2 {\rm (cf.[6])}} The pair $(V_1,V_2)$ of linear spaces
is called {\it an anti--Jordan pair\/} iff there are defined two mappings
$m_1:V_2\otimes\bigwedge^2V_1\mapsto V_1$ and
$m_2:V_1\otimes\bigwedge^2V_2\mapsto V_2$ such that the mappings $(X,Y)\mapsto
[X,Y]_A=m_1(A,X,Y)$ ($X,Y\in V_1$, $A\in V_2$) and $(A,B)\mapsto
[A,B]_X=m_2(X,A,B)$ ($A,B\in V_2$, $X\in V_1$) are compatible to each other
in the following manner
$$\align
[X,Y]_{[A,B]_Z}=&[[X,Z]_A,Y]_B+[[Z,Y]_A,X]_B-[[X,Y]_B,Z]_A\endalign$$ and
$$\align
[A,B]_{[X,Y]_C}=&[[A,C]_X,B]_Y+[[C,B]_X,A]_Y-[[A,B]_Y,C]_X\endalign$$
($X,Y,Z\in V_1$, $A,B,C\in V_2$) hold.
\enddefinition

It can be easily verified that isotopic pairs are always anti--Jordan pairs
(to obtain it one should use the Jacobi identity linearized by subscript
parameters), and that the anti--Jordan pairs with a multiplication, obeying
Jacobi identity if a subscript parameter is fixed in any way, are just the
isotopic pairs.  So the isotopic pairs may be considered as a particular case
of the anti--Jordan pairs. Note that there exist examples of anti--Jordan
pairs, which are not isotopic ones [5].

Anti--Jordan pairs are closely related to the (polarized) anti--Lie triple
systems and Lie superalgebras [5] (cf. also [7]). Namely,

\definition{Definition 3} The ternary algebra $V$ with product $[xyz]$ is
called {\it an anti--Lie triple system\/} if
$$\align
&(1)\quad [xyz]=[xzy],\\
&(2)\quad [xyz]+[zxy]+[yzx]=0,\\
&(3)\quad [[xyz]uv]=[[xuv]yz]+[x[yvu]z]+[xy[zuv]].
\endalign$$
An anti--Lie triple system $V$ is {\it polarized\/} iff $V=V_1\oplus
V_2$ and $[xyz]=0$ for $y,z\in V_1$ or $y,z\in V_2$.
\enddefinition

If $V$ is an anti--Lie triple system let's put $R_{yz}\in\End(V):
R_{yz}x=[xyz]$. The operators $R_{yz}$ are closed under commutators so that
$\frak g_0(V)=\spn(R_{yz}; y,z\in V)$ is a Lie algebra. The space $\frak
g_0(V)\oplus V$ possesses a natural structure of a Lie superalgebra [8] with
the even part $\frak g_0(V)$ and the odd part $V$ [5]. It will be denoted by
$\frak g(V)$. Polarized anti--Lie triple systems $V=V_1\oplus V_2$ produce
polarized Lie superalgebras $\frak g(V)=\frak g_0(V)\oplus(V_1\oplus V_2)$
such that $[V_i,V_i]_+=0$, $[\frak g(V), V_i]_-\subseteq V_i$ (it should be
marked that there is sometimes asserted that $V_2\simeq V_1^*$ as $\frak
g_0(V)$--modules, however, we shall not do it in general).

An arbitrary anti--Jordan pair (so an isotopic pair, in particular) has a
structure of a polarized anti--Lie triple system. Namely, one should put
$[xyz]=[z,x]_y$ (iff $z$ belongs to the same space $V_i$ as $x$) and $[y,x]_z$
(iff $y$ belongs to the same space $V_i$ as $x$).

\remark{Remark 3} Let's summarize the relations between the concepts of
"isotopic pair", "anti--Jordan pair", "polarized anti--Lie triple system" and
"polarized Lie superalgebra" once more.
\roster
\item"(1)" Each isotopic pair is an anti--Jordan pair, though there exist
anti--Jordan pairs, which are not isotopic. It means that isotopic pairs form
a proper subclass of the class of anti--Jordan pairs.
\item"(2)" Categories of anti--Jordan pairs and polarized anti--Lie triple
systems are equivalent. It means that each anti--Jordan pair defines a
polarized anti--Lie triple system and vice versa.
\item"(3)" Categories of polarized anti--Lie triple systems and polarized Lie
superalgebras are equivalent. On the other hand polarized anti--Lie triple
systems and polarized Lie superalgebras are particular cases of anti--Lie
triple systems and Lie superalgebras respectively and the marked equivalency of
categories is a particular case of the equivalency of categories of anti--Lie
triple systems and Lie superalgebras.
\item"(4)" As a consequence of (2) and (3) categories of anti--Jordan pairs
and polarized Lie superalgebras are equivalent.
\item"(5)" As a consequence of (1) and (4) each isotopic pair defines a
polarized Lie superalgebra but not vice versa.
\endroster\endremark

An illustrative example to the construction of a Lie superalgebra by an
isotopic pair is convenient. {\it Example\/}: let $H_1$ and $H_2$ be two
linear spaces, $(\Hom(H_1,H_2);\mathbreak\Hom(H_2,H_1))$ is an isotopic pair,
the corresponding Lie superalgebra is isomorphic to $\operatorname{\frak
g\frak l}(n|m)$, $n=\dim H_1$, $m=\dim H_2$.

Note that the isocommutators in an isotopic pair $(V_1,V_2)$
define families of Poisson brackets $\{\cdot,\cdot\}_A$ and
$\{\cdot,\cdot\}_X$ ($A\in V_2$, $X\in V_1$) in the spaces $S^{\cdot}(V_1)$
and $S^{\cdot}(V_2)$, respectively.

\definition{Definition 4 {\rm (cf.[3])}} Let's consider two elements $\Cal
H_1$ and $\Cal H_2$ ({\it "hamiltonians"\/}) in $S^{\cdot}(V_1)$ and
$S^{\cdot}(V_2)$, respectively. The equations
$$\dot X_t=\{\Cal H_1,X_t\}_{A_t},\qquad \dot A_t=\{\Cal H_2,A_t\}_{X_t},$$
where $X_t\in V_1$ and $A_t\in V_2$ are called {\it the (nonlinear) dynamical
equations associated with the isotopic pair $(V_1,V_2)$ and "hamiltonians"
$\Cal H_1$ and $\Cal H_2$\/} (it should be marked that "hamiltonians" are not
even integrals of motion in a general situation).
\enddefinition

It should be mentioned that the dynamical equations associated with isotopic
pairs are a particular case of such equations associated with general I--pairs
[4].

\head III. ISOTOPIC PAIR OF NONCANONICALLY COUPLED OSCILLATORS:
ALGEBRAIC ASPECTS
\endhead

Let's now consider the isotopic pairs of noncanonically coupled oscillators
[3,6]. The space $V_1$ is spanned by the elements $p$, $q$ and $r$ and the
space $V_2$ is spanned by the elements $a$, $b$ and $c$. The isocommutators
have the form

\

\centerline{
$\aligned
[p,q]_a&=2\varepsilon_1 q\\
[p,r]_a&=\varepsilon_2 r\\
[q,r]_a&=0
\endaligned$
$\quad$
$\aligned
[p,q]_b&=2\varepsilon_1 p\\
[p,r]_b&=0\\
[q,r]_b&=-\varepsilon_2 r
\endaligned$
$\quad$
$\aligned
[p,q]_c=\varepsilon_3 r\\
[p,r]_c=0\\
[q,r]_c=0
\endaligned$}

\

\

\centerline{
$\aligned
[a,b]_p&=2\tilde\varepsilon_1 b\\
[a,c]_p&=\tilde\varepsilon_2 c\\
[b,c]_p&=0
\endaligned$
$\quad$
$\aligned
[a,b]_q&=2\tilde\varepsilon_1 a\\
[a,c]_q&=0\\
[b,c]_q&=-\tilde\varepsilon_2 c
\endaligned$
$\quad$
$\aligned
[a,b]_r=\tilde\varepsilon_3 c\\
[b,c]_r=0\\
[a,c]_r=0
\endaligned$}

\

where $$\left\{\aligned &\varepsilon_1+\tilde\varepsilon_1=0\\
&\varepsilon_2-\tilde\varepsilon_2=\varepsilon_1-\tilde\varepsilon_1\\
&\varepsilon_3\tilde\varepsilon_3-\varepsilon_2\tilde\varepsilon_2=0\endaligned
\right.$$

The corresponding Lie algebra $\frak g_0(V_1\oplus V_2)$ is spanned (for
generic $\varepsilon_i$, $\tilde\varepsilon_i$) by 6 operators
$R_{p,a}$, $R_{p,b}$, $R_{q,a}$, $R_{q,b}$,
$R_{r,b}=\frac{\varepsilon_2}{\varepsilon_3}R_{p,c}$,
$R_{r,a}=\frac{\varepsilon_2}{\varepsilon_3}R_{q,c}$,
which have the form
$$\align
R_{p,a}=\left(\matrix 2\varepsilon_1 & 0 & 0 \\ 0 & 0 & 0 \\ 0 & 0 &
\varepsilon_2\endmatrix\right),& \quad
R_{p,b}=\left(\matrix 0 & 0 & 0 \\ 2\varepsilon_1 & 0 & 0 \\ 0 & 0 &
0\endmatrix\right),\\
R_{q,a}=\left(\matrix 0 & -2\varepsilon_1 & 0 \\ 0 & 0 & 0 \\ 0 & 0 &
0\endmatrix\right),& \quad
R_{q,b}=\left(\matrix 0 & 0 & 0 \\ 0 & -2\varepsilon_1 & 0 \\ 0 & 0 &
-\varepsilon_2\endmatrix\right),\\
R_{p,c}=\left(\matrix 0 & 0 & 0 \\ 0 & 0 & 0 \\ \varepsilon_3 & 0 &
0\endmatrix\right),& \quad
R_{q,c}=\left(\matrix 0 & 0 & 0 \\ 0 & 0 & 0 \\ 0 & -\varepsilon_3 &
0\endmatrix\right)
\endalign $$
in the basis $(q,p,r)$ and the form
$$\align
R_{p,a}=\left(\matrix 0 & 0 & 0 \\ 0 & 2\tilde\varepsilon_1 & 0 \\ 0 & 0 &
\tilde\varepsilon_2\endmatrix\right),& \quad
R_{p,b}=\left(\matrix 0 & 0 & 0 \\ -2\tilde\varepsilon_1 & 0 & 0 \\ 0 & 0 &
0\endmatrix\right),\\
R_{q,a}=\left(\matrix 0 & 2\tilde\varepsilon_1 & 0 \\ 0 & 0 & 0 \\ 0 & 0 &
0\endmatrix\right),& \quad
R_{q,b}=\left(\matrix -2\tilde\varepsilon_1 & 0 & 0 \\ 0 & 0 & 0 \\ 0 & 0 &
-\tilde\varepsilon_2\endmatrix\right),\\
R_{p,c}=\left(\matrix 0 & 0 & 0 \\ 0 & 0 & 0 \\ -\tilde\varepsilon_2 & 0 &
0\endmatrix\right),& \quad
R_{q,c}=\left(\matrix 0 & 0 & 0 \\ 0 & 0 & 0 \\ 0 & \tilde\varepsilon_2 &
0\endmatrix\right)
\endalign $$
in the basis $(a,b,c)$.

The Lie superalgebra $\frak g(V_1\oplus V_2)$ has a (super)dimension $(6|6)$
and is generated by $R_{p,a}$, $R_{p,b}$, $R_{q,a}$, $R_{q,b}$, $R_{p,c}$,
$R_{q,c}$, $p$, $q$, $r$, $a$, $b$, $c$ with (super)commutation relations
$$\align
[q,p]_+=[q,r]_+=[p,r]_+=[a,b]_+&=[a,c]_+=[b,c]_+=[r,c]_+=0,\\
[p,a]_+=R_{p,a},\ [q,a]_+&=R_{q,a},\ [p,b]_+=R_{p,b},\\ [q,b]_+=R_{q,b},\
[p,c]_+&=R_{p,c},\ [q,c]_+=R_{q,c},\\
[r,a]_+=\tfrac{\varepsilon_2}{\varepsilon_3}R_{q,c},&\
[r,b]_+=\tfrac{\varepsilon_2}{\varepsilon_3}R_{p,c};\\
&\\ \allowdisplaybreak
[R_{p,a},q]_-=2\varepsilon_1q,\ [R_{p,a},p]_-&=0,\
[R_{p,a},r]_-=\varepsilon_2r,\\
[R_{q,a},q]_-=0,\ [R_{q,a},p]_-&=-2\varepsilon_1q,\ [R_{q,a},r]_-=0,\\
[R_{p,b},q]_-=2\varepsilon_1p,\ [R_{p,b},p]_-&=0,\ [R_{p,b},r]_-=0,\\
[R_{q,b},q]_-=0,\ [R_{q,b},p]_-&=-2\varepsilon_1p,\
[R_{q,b},r]_-=-\varepsilon_2r,\\
[R_{p,c},q]_-=\varepsilon_3r,\ [R_{p,c},p]_-&=0,\ [R_{p,c},r]_-=0,\\
[R_{q,c},q]_-=0,\ [R_{q,c},p]_-&=-\varepsilon_3r,\ [R_{q,c},r]_-=0,\\
[R_{p,a},a]_-=0,\ [R_{p,a},b]_-&=2\tilde\varepsilon_1b,\
[R_{p,a},c]_-=\tilde\varepsilon_2c,\\
[R_{q,a},a]_-=0,\ [R_{q,a},b]_-&=2\tilde\varepsilon_1a,\ [R_{q,a},c]_-=0,\\
[R_{p,b},a]_-=-2\tilde\varepsilon_1b,\ [R_{p,b},b]_-&=0,\ [R_{p,b},c]_-=0,\\
[R_{q,b},a]_-=-2\tilde\varepsilon_1a,\ [R_{q,b},b]_-&=0,\
[R_{q,b},c]_-=-\tilde\varepsilon_2c,\\
[R_{p,c},a]_-=-\tilde\varepsilon_2c,\ [R_{p,c},b]_-&=0,\ [R_{p,c},c]_-=0,\\
[R_{q,c},a]_-=0,\ [R_{q,c},b]_-&=\tilde\varepsilon_2c,\ [R_{q,c},c]_-=0;\\
&\\ \allowdisplaybreak
[R_{p,a},R_{p,b}]_-=-2\varepsilon_1R_{p,b},\
[R_{p,a},R_{q,a}]_-&=2\varepsilon_1R_{q,a},\ [R_{p,a},R_{p,b}]_-=0,\\
[R_{p,a},R_{p,c}]_-=\tilde\varepsilon_2R_{p,c},\
[R_{p,a},R_{q,c}]_-&=\varepsilon_2R_{q,c},\
[R_{p,b},R_{q,a}]_-=2\varepsilon_1(R_{q,b}+R_{p,a}),\\
[R_{p,b},R_{q,b}]_-=2\varepsilon_1R_{p,b},\ [R_{p,b},R_{p,c}]_-&=0,\
[R_{p,b},R_{q,c}]_-=2\varepsilon_1R_{p,c},\\
[R_{q,a},R_{q,b}]_-=-2\varepsilon_1R_{q,a},\
[R_{q,a},R_{p,c}]_-&=-2\varepsilon_1R_{p,c},\
[R_{q,a},R_{q,c}]_-=0,\\
[R_{q,b},R_{p,c}]_-=-\varepsilon_2R_{p,c},\
[R_{q,b},R_{q,c}]_-&=-\tilde\varepsilon_2R_{q,c},\ [R_{p,c},R_{q,c}]_-=0.
\endalign$$

The even part of the Lie superalgebra $\frak g(V_1\oplus V_2)$ is isomorphic
to the semidirect sum of $\operatorname{\frak g\frak l}(2,\Bbb C)$ and $\Bbb
C^2$. On the other hand $\frak g(V_1\oplus V_2)$ may be considered as a
semidirect product of the Lie superalgebra $\operatorname{\frak s\frak
l}(2|1,\Bbb C)$ generated by $R_{p,a}$, $R_{p,b}$, $R_{q,a}$, $R_{q,b}$, $p$,
$q$,
$a$, $b$ and the $(2|2)$--dimensional vector superspace $V^{2|2}$ generated by
$R_{p,c}$, $R_{q,c}$, $r$, $c$.

\head IV. ISOTOPIC PAIR OF NONCANONICALLY COUPLED OSCILLATORS:
CLASSICAL DYNAMICS
\endhead

The dynamical equations with "hamiltonians" $\Cal H_1=P^2+Q^2$ and $\Cal
H_2=A^2+B^2$ have the form
$$\left\{\aligned
\dot P=&-4\varepsilon_1(Q^2A+PQB)-2\varepsilon_3RQC\\
\dot Q=&4\varepsilon_1(PQA+P^2B)+2\varepsilon_3RPC\\
\dot R=&2\varepsilon_2(PRA-QRB)\endaligned\right.$$
$$\left\{\aligned
\dot A=&-4\tilde\varepsilon_1(B^2P+ABQ)-2\tilde\varepsilon_3CBR\\
\dot B=&4\tilde\varepsilon_1(ABP+A^2Q)+2\tilde\varepsilon_3CAR\\
\dot C=&2\tilde\varepsilon_2(ACP-BCQ)\endaligned\right.$$

Note that "hamiltonians" $\Cal H_1=\Cal I_1^2$ and $\Cal H_2=\Cal I_2^2$
are integrals of motion here, so it is rather convenient to put $P=\Cal
I_1\cos\varphi$, $Q=\Cal I_1\sin\varphi$, $A=\Cal I_2\cos\psi$, $B=\Cal
I_2\sin\psi$. Then
$$\left\{\aligned
\dot\varphi=&-2\varepsilon_3RC-4\varepsilon_1\Cal I_1\Cal
I_2\sin(\varphi+\psi)\\
\dot\psi=&-2\tilde\varepsilon_3RC-4\tilde\varepsilon_1\Cal I_1\Cal
I_2\sin(\varphi+\psi)\endaligned\right.$$
$$\left\{\aligned
\dot R=&2\varepsilon_2\cos(\varphi+\psi)R\\
\dot C=&2\tilde\varepsilon_2\cos(\varphi+\psi)C
\endaligned\right.$$
Let's introduce $\vartheta=\varphi+\psi$, $\chi=\varepsilon_3\psi-
\tilde\varepsilon_3\varphi$ and mark that $\varepsilon_1+
\tilde\varepsilon_1=0$, then
$$\left\{\aligned
\dot\vartheta=&-2(\varepsilon_3+\tilde\varepsilon_3)RC\\
\dot\chi=&4\varepsilon_1\Cal I_1\Cal I_2(\varepsilon_3-\tilde\varepsilon_3)
\sin\vartheta\endaligned\right.$$
Also
$$(RC)^{\cdot}=2(\varepsilon_2+\tilde\varepsilon_2)\cos\vartheta(RC),$$
therefore,
$$(RC)'_{\vartheta}=-\frac{\varepsilon_2+\tilde\varepsilon_2}
{\varepsilon_3+\tilde\varepsilon_3}\cos\vartheta$$
and
$$RC=\Cal L-\frac{\varepsilon_2+\tilde\varepsilon_2}{\varepsilon_3+\tilde
\varepsilon_3}\Cal I_1\Cal I_2\sin\vartheta,$$
whereas
$$\dot\vartheta=-2\Cal L(\varepsilon_3+\tilde\varepsilon_3)-2\Cal I_1\Cal I_2
(\varepsilon_2+\tilde\varepsilon_2)\sin\vartheta.$$
Here $\Cal L=RC+\frac{\varepsilon_2+\tilde\varepsilon_2}{\varepsilon_3+\tilde
\varepsilon_3}(QA+PB)$ is an integral of motion.
Note that
$$(R^{\tilde\varepsilon_2}C^{-\varepsilon_2})^{\cdot}=0$$
so it is convenient to put
$$\Lambda=R^{\frac{\tilde\varepsilon_2}{\varepsilon_2+\tilde\varepsilon_2}}
C^{\frac{\varepsilon_2}{\varepsilon_2+\tilde\varepsilon_2}}.$$
Then
$$\left\{\aligned
R=&\Lambda(\Cal L-\frac{\varepsilon_2+\tilde\varepsilon_2}{\varepsilon_3+\tilde
\varepsilon_3}\Cal I_1\Cal I_2\sin\vartheta)^{\frac{\varepsilon_2}
{\varepsilon_2+\tilde\varepsilon_2}}\\
C=&\frac1{\Lambda}(\Cal
L-\frac{\varepsilon_2+\tilde\varepsilon_2}{\varepsilon_3+
\tilde\varepsilon_3}\Cal I_1\Cal I_2\sin\vartheta)^{\frac{\tilde\varepsilon_2}
{\varepsilon_2+\tilde\varepsilon_2}}\endaligned\right.$$
$\Cal I_1$, $\Cal I_2$, $\Cal L$ and $\Lambda$ form a complete set of
integrals of motion for generic values of $\varepsilon_i$,
$\tilde\varepsilon_i$.

Let's also denote
$$\align\xi=(\varepsilon_2&+\tilde\varepsilon_2)\chi+
2\varepsilon_1(\varepsilon_3-\tilde\varepsilon_3)\vartheta\\
&=[(\varepsilon_2+\tilde\varepsilon_2)\varepsilon_3+
2\varepsilon_1(\varepsilon_3-\tilde\varepsilon_3)]\psi-[(\varepsilon_2+
\tilde\varepsilon_2)\tilde\varepsilon_3+2\tilde\varepsilon_1(\varepsilon_3-
\tilde\varepsilon_3)]\varphi,\endalign$$
then
$$\xi=4\Cal L(\tilde\varepsilon^2_3-\varepsilon^2_3)\varepsilon_1t+\xi_0.$$

The obtained results certainly generalize results of [3]. It is also possible
to consider a hybrid coupling for the isotopic pairs of noncanonically coupled
oscillators (cf.[3]).

\head V. ISOTOPIC PAIR OF NONCANONICALLY COUPLED OSCILLATORS: REPRESENTATION
THEORY AND QUANTUM DYNAMICS
\endhead

\definition{Definition 5}
{\it A representation of the isotopic pair
$(V_1,V_2)$ in the linear space $W$\/} is a pair $(T_1,T_2)$ of mappings
$T_i:V_i\mapsto\End(W)$ such that
$$\align
T_1([X,Y]_A)=&T_1(X)T_2(A)T_1(Y)-T_1(Y)T_2(A)T_1(X),\\
T_2([A,B]_X)=&T_2(A)T_1(X)T_2(B)-T_2(B)T_1(X)T_2(A),
\endalign$$
where $X,Y\in V_1$, $A,B\in V_2$ [6]. A representation of the isotopic pair
$(V_1,V_2)$ in the linear space $W$ is called {\it nilpotent\/} if
$$\left\{\aligned (\forall X_i\in V_1) &\ T_1(X_1)T_1(X_2)=0,\\
(\forall A_i\in V_2) &\ T_2(A_1)T_2(A_2)=0.
\endaligned\right.
$$
A representation of the isotopic pair $(V_1,V_2)$ in
the linear space $W$ is called {\it split\/} [6] iff $W=W_1\oplus W_2$ and
$$\left\{\aligned
(\forall X\in V_1) \left. T_1(X)\right|_{W_2}=0,\ T_1(X):W_1\mapsto W_2,\\
(\forall A\in V_2) \left. T_2(A)\right|_{W_1}=0,\ T_2(A):W_2\mapsto W_1.
\endaligned\right.$$
Otherwords, operators $T(X)$ and $T(A)$ have the form $\left(\matrix 0 & 0
\\ * & 0 \endmatrix\right)$ and $\left(\matrix 0 & * \\ 0 & 0 \endmatrix
\right)$, respectively. Each split representation is nilpotent.
\enddefinition

Not that a nilpotent representation (but not an arbitrary one) of an isotopic
pair $(V_1,V_2)$ defines a representation $T$ of the corresponding anti--Lie
triple system and Lie superalgebra $\frak g(V_1\oplus V_2)$ (or its central
extension $\hat{\frak g}(V_1\oplus V_2)$). It should be stressed once more
that a representation of the polarized Lie superalgebra, constructed by an
isotopic pair, may not define a representation of the least. If the
representation of the isotopic pair $(V_1,V_2)$ is split then the
representation of the Lie superalgebra $\frak g(V_1\oplus V_2)$ always have a
special "polarized" form: $W=W_1\oplus W_2$, $T(V_1):W_1\mapsto W_2$,
$T(V_2):W_2\mapsto W_1$, $\frak g_0(V_1\oplus V_2):W_i\mapsto W_i$.

Note that each representation $(T_1,T_2)$ of the isotopic pair $(V_1,V_2)$ in
the space $W$ defines a split representation $(T^s_1,T^s_2)$ of the same pair
in the space $W_1\oplus W_2$ ($W_i\simeq W$):
$$(\forall X\in V_1)\ T^s_1(X)=\left(\matrix 0 & 0\\ T_1(X) & 0
\endmatrix\right),\quad (\forall A\in V_2)\ T^s_2(A)=\left(\matrix 0 & T_2(A)\\
0 & 0 \endmatrix\right).$$

To construct a split representation of an isotopic pair
$(V_1,V_2)$ one may start from two arbitrary $\frak g_0(V_1\oplus
V_2)$--modules $W_1$ and $W_2$, to consider suitable $\frak g_0(V_1\oplus
V_2)$--tensor operators [9] from $W_1$ to $W_2$ and vice versa as candidates
for $T(V_1)$ and $T(V_2)$, respectively, and then to check the validity of
anticommutation relations between the tensor operators. In such approach
elements of the isotopic pair $(V_1,V_2)$ are realized as hidden symmetries
with respect to $\frak g_0(V_1\oplus V_2)$ (cf.[10,11]).

Let's now describe the quantum dynamics of noncanonically coupled
oscillators.

Mark that the classical dynamics of noncanonically coupled oscillators was not
somehow naturally related to Lie superalgebras, the situation in the quantum
case is rather different. Namely, the formal quantum dynamical equations have
the form
$$\left\{\aligned
\tfrac{d}{dt}\P&=-2\varepsilon_1(\P\B\Q+\Q\B\P+2\Q\A\Q)-
\varepsilon_3(\R\C\Q+\Q\C\R)\\
\tfrac{d}{dt}\Q&=2\varepsilon_1(\P\A\Q+\Q\A\P+2\P\B\P)+
\varepsilon_3(\R\C\P+\P\C\R)\\
\tfrac{d}{dt}\R&=\varepsilon_2(\P\A\R+\R\A\P-\Q\B\R-\R\B\Q)
\endaligned\right.$$

$$\left\{\aligned
\tfrac{d}{dt}\A&=-2\tilde\varepsilon_1(\A\Q\B+\B\Q\A+2\B\P\B)-
\tilde\varepsilon_3(\C\R\B+\B\R\C)\\
\tfrac{d}{dt}\B&=2\tilde\varepsilon_1(\A\P\B+\B\P\A+2\A\Q\A)+
\tilde\varepsilon_3(\C\R\A+\A\R\C)\\
\tfrac{d}{dt}\C&=\tilde\varepsilon_2(\A\P\C+\C\P\A-\B\Q\C-\C\Q\B)
\endaligned\right.$$

The dynamics is considered in arbitrary representation of the isotopic pair of
noncanonically coupled oscillators. Let's consider such dynamics in the
corresponding split representation. First of all renormalize $c$ and $r$ so
that $R_{p,c}=R_{b,r}$ and $R_{q,c}=R_{a,r}$. Then the following proposition
holds.

\define\hi{\operatorname{hidden}}
\proclaim{Proposition 1} Equations of quantum dynamics of noncanonically
coupled oscillators are a reduction of formal super Heisenberg equations
$$\tfrac{d}{dt}\hat F_t=[\hat H_{\hi},\hat F_t]$$
in $\Cal U(\frak g(V_1\oplus V_2))$ with quadratic quantum hamiltonian
$$\hat H_{\hi}=\hat R_{q,a}^2+\hat R_{p,b}^2+\hat R_{q,b}^2+\hat R_{p,a}^2+
\hat R_{p,c}^2+\hat R_{q,c}^2$$
\endproclaim

So quantum dynamics of noncanonically coupled oscillators admits a hidden
super--hamiltonian formulation in terms of Lie superalgebra $\frak g(V_1\oplus
V_2)$.

\proclaim{Corollary} The quantum dynamics preserves the initial operator
relations:
$$\aligned &\P\A\Q-\Q\A\P=2\varepsilon_1\Q,\ \P\A\R-\R\A\P=\varepsilon_2\R,\
\Q\A\R-\R\A\Q=0,\\
&\P\B\Q-\Q\B\P=2\varepsilon_1\P,\ \P\B\R-\R\B\P=0,\
\Q\B\R-\R\B\Q=-2\varepsilon_2\R,\\
& \P\C\Q-\Q\C\P=\varepsilon_3\R,\ \P\C\R-\R\C\P=0,\ \Q\C\R-\R\C\Q=0,\\
&\\
& \A\P\B-\B\P\A=2\tilde\varepsilon_1\B,\ \A\P\C-\C\P\A=\tilde\varepsilon_2\C,\
\B\P\C-\C\P\B=0,\\
& \A\Q\B-\B\Q\A=2\tilde\varepsilon_1\A,\ \A\Q\C-\C\Q\A=0,\
\B\Q\C-\C\Q\B=-\tilde\varepsilon_2\C,\\
& \A\R\B-\B\R\A=\tilde\varepsilon_3\C,\ \B\R\C-\C\R\B=0,\ \A\R\C-\C\R\A=0.
\endaligned
$$
\endproclaim

\head VI. REMARKS \endhead

\remark{Remark 4} Note that the dynamics of the classical Poisson brackets is
not defined {\it as such}, because the subscript parameters are considered to
be only linear functions. Hence, to define their conserved dynamics it is
necessary to extend them correctly to arbitrary subscript parameters.
\endremark

\remark{Remark 5} It is interesting to consider a quantum version of hybrid
couplings [3].
\endremark

\remark{Remark 6} Though the representation theory of isotopic pairs may be
imbed into the representation theory of (polarized) Lie superalgebras,
sometimes it is rather reasonable to avoid it in view of the dimension growing
(f.e. for isotopic pairs of symmetric and skew--symmetric matrices [3,5] or
for geometric isotopic pairs [5]).
\endremark

\head VII. CONCLUSIONS \endhead

Thus, the classical and quantum dynamics of noncanonically coupled oscillators
is investigated. The crucial role of Lie superalgebras is explicated. It is
shown that quantum dynamics preserves the initial operator relations.

\head ACKNOWLEDGEMENTS
\endhead

I thank Laboratoire de Physique Th\'eorique, ENS, Paris for beautiful
atmosphere and perfect conditions as well as Max--Plank Institut, Bonn,
Germany for a support during my visit there in September 1994, where the work
was finished.

\head APPENDIX A: LIE $\frak g$--BUNCHES AND ISOTOPIC PAIRS
\endhead

\definition{Definition 6}

{\bf A [6].} Let $\frak g$ be a Lie algebra. A {\it Lie $\frak g$--bunch\/} is
a $\frak g$--module $V$ such that there exists a $\frak g$--equivariant
mapping $\frak g\otimes\bigwedge^2(V)\mapsto V$, which defines a structure of
Lie algebra in $V$ when the first argument is fixed in an arbitrary way; we
shall denote this mapping by $[\cdot,\cdot]_A$, $A\in\frak g$.

{\bf B.} A Lie $\frak g$--bunch $V$ is called {\it complete\/} if the Lie
brackets defined by elements of $\frak g$ are closed under the $\lozenge$
operation, i.e.
$$\aligned \forall A,B\in\frak g\ \forall X\in V\ \exists
C=A\underset X\to\lozenge B\in\frak g:&\\
[X,Y]_C=\tfrac12([[X,Z]_A,Y]_B+&[[X,Y]_A,Z]_B+[[Z,Y]_A,X]_B-
(\alpha\longleftrightarrow\beta)).\endaligned$$
\enddefinition

Numerous examples of Lie $\frak g$--bunches are collected in [9].

\proclaim{Proposition 2} Each complete Lie $\frak g$--bunch $V$ may be enlarged
to an isotopic pair $(V\oplus\Bbb C,\frak g)$ with isocommutators
$$\aligned
[(X,\lambda\boldkey 1),(Y,\mu\boldkey 1)]_A&=([X,Y]_A+\lambda A(Y)-\mu
A(X),0)\\
[A,B]_{(X,\lambda\boldkey 1)}&=A\underset X\to\lozenge B+\lambda[A,B]
\endaligned$$
where $X,Y\in V$, $A,B\in\frak g$; $A(X)$ is an image of an element $X\in V$
under the action of an element $A\in\frak g$, $[A,B]$ is the
commutator in $\frak g$ and $[X,Y]_A$ is the isocommutator in $V$.
\endproclaim

The proof of the proposition is straightforward.

\head APPENDIX B: ISOREPRESENTATIONS OF LIE ALGEBRAS VIA LIE SUPERALGEBRAS
\endhead

\definition{Definition 7} Let $\frak g$ be a Lie algebra. {\it An
isorepresentation\/} ({\it $Q$--isorepresentation\/}) of $\frak g$ in the
linear space $V$ is a pair $(T,Q)$, $T\in Hom(\frak g,\End(V))$, $Q\in\End(V)$
such that $\forall X,Y\in\frak g$ $T(X)QT(Y)-T(Y)QT(X)=T([X,Y])$.
\enddefinition

Note that a Lie algebra $\frak g$ maybe considered as a complete Lie $\Bbb
C$--bunch. The corresponding isotopic pair $(\Bbb C,\frak g)$ will be denoted
by $\Cal I(\frak g)$. The isorepresentations of $\frak g$ are just the
representations of $\Cal I(\frak g)$. The Lie superalgebra constructed from
$\Cal I(\frak g)$ is a standard Lie superalgebra associated with Lie algebra
$\frak g$ with even part isomorphic to $\frak g$, which acts in the odd part
as in $\ad_{\frak g}\oplus\boldkey 1_{\frak g}$ ($\ad_{\frak g}$ is the
adjoint $\frak g$--module and $\boldkey 1_{\frak g}$ is the trivial
one--dimensional $\frak g$--module). The mapping from the symmetric square of
the odd part to the even part is natural: $\ad_{\frak g}\otimes\boldkey
1_{\frak g}\mapsto\frak g$. It immediately gives a description of irreducible
representations of $\frak g$. First of all, let's consider a crucial example:

\define\GL{\operatorname{GL}}
\remark{Example} Let $(T_Q,Q)$ be an isorepresentation of
$\frak g$ in the space $V$. Then $V$ admits two natural representation
$T^{\pm}$ of the Lie algebra $\frak g$: $\forall X\in\frak g$ $T^+(X)=QT_Q(X)$
and $T^-(X)=T_Q(X)Q$, which are equivalent if $Q$ is invertible.\newline
Let $T$ be a representation of $\frak g$ in the space $V$ and $Q\in\GL(V)$,
then one may construct two equivalent isorepresentations $(T^{\pm}_Q,Q)$ of
$\frak g$ in $V$: $\forall X\in\frak g$ $T^+_Q(X)=Q^{-1}T(X)$ and
$T^-_Q(X)=T(X)Q^{-1}$.
\endremark

Second, one may restrict himself by split isorepresentations (i.e. split
representations of $\Cal I(\frak g)$).

\proclaim{Proposition 3} Each finite dimensional irreducible split
isorepresentation of $\frak g$ may be realized in the linear space
$V=V_1\oplus V_2$ ($V_1\simeq V_2\simeq V_{\alpha}$, where $V_{\alpha}$
carries an irreducible representation $T_{\alpha}$ of $\frak g$) and
$Q$ is an isomorphism of $V_1$ onto $V_2$.
\endproclaim

\remark{Example {\rm (Two--dimensional irreducible isorepresentation of $\Bbb
C$, cf.[9])}} The operator $Q$ have the form $\left(\matrix 0 & 0 \\ 1 &
0\endmatrix\right)$, the generator of $\Bbb R$ is represented by $\left(\matrix
0 &
1 \\ 0 & 0 \endmatrix\right)$.
\endremark

This example may be straightforwardly generalized on an arbitrary Lie algebra
$\frak g$. Namely, let $V=\ad_{\frak g}\oplus\ad_{\frak g}$, the operator $Q$
has the form $\left(\matrix 0 & 0 \\ E & 0 \endmatrix\right)$, where $E$ is
a unit matrix and an element $X$ of $\frak g$ is represented by $\left(\matrix
0 & \ad_{\frak g}(X) \\ 0 & 0\endmatrix\right)$.

\Refs
\roster
\item"[1]" Arnold V I 1976 {\it Mathematical methods of classical mechanics},
Springer--Verlag;\newline
Dubrovin B A, Novikov S P and Fomenko A T 1988 {\it Modern geometry --- methods
and applications}, Springer--Verlag;\newline
Karasev M V and Maslov V P 1993 {\it Nonlinear Poisson brackets. Geometry and
quantization}, Amer.~Math.~Soc., RI.
\item"[2]" Alekseev A Yu and Faddeev L D 1991 {\it Commun.~Math.~Phys.\/} {\bf
141} 413-422.
\item"[3]" Juriev D 1994 Topics in nonhamiltonian interaction of hamiltonian
dynamic systems,\linebreak {\it mp\_arc/94-136}.
\item"[4]" Juriev D 1994 On the nonhamiltonian interaction of two rotators,
{\it dg-ga/9409004}.
\item"[5]" Faulkner J R 1973 {\it J.~Algebra\/} {\bf 26} 1-9;\newline
Faulkner J R and Ferrar J C 1980 {\it Commun.~Alg.\/} {\bf 8}
993-1013.
\item"[6]" Juriev D 1994 Topics in isotopic pairs and their representations,
{\it mp\_arc/94-267}.
\item"[7]" Okubo S 1994 {\it J.~Math.~Phys.\/} {\bf 35} 2785-2803.
\item"[8]" Kac V G 1977 {\it Commun.~Math.~Phys.\/} {\bf 53} 31-64; 1977
{\it Adv.~Math.\/} {\bf 26} 8-96; 1978 {\it Lect.~Notes Math.\/} {\bf 676}
597-626;\newline
Scheunert M 1979 {\it The theory of Lie superalgebras}.
Springer--Verlag;\newline
Manin Yu I 1988 {\it Gauge field theory and complex geometry}.
Springer--Verlag.
\item"[9]" Biedenharn L C and Louck J D 1981 {\it Angular momentum in quantum
mechanics. Theory and applications}. Addison Wesley Publ; {\it The
Racah--Wigner algebra in quantum theory}. Addison Wesley Publ.
\item"[10]" Juriev D 1994 {\it J.~Math.~Phys.\/} {\bf 35\rm:9}.
\item"[11]" Juriev D 1994 Topics in hidden symmetries, {\it hep-th/9405050}.
\endroster
\endRefs
\enddocument